УДК 517.956.6;517.44         DOI: 10.20998/2411-0558.2017.50.13

*Т.Г. ЭРГАШЕВ*, канд.физ.-мат.наук, доц., Ташкентский институт инженеров ирригации и механизации сельского хозяйства, г. Ташкент, Узбекистан

## ФОРМУЛА ОБРАЩЕНИЯ ИНТЕГРАЛЬНОГО УРАВНЕНИЯ ВОЛЬТЕРРА С ФУНКЦИЕЙ ГУМБЕРТА В ЯДРЕ И ЕЁ ПРИЛОЖЕНИЯ К РЕШЕНИЮ КРАЕВЫХ ЗАДАЧ

Многие задачи прикладной математики сводятся к решению интегральных уравнений со специальными функциями в ядрах, поэтому формулы обращения таких уравнений играют важную роль при решении различных задач. В работе введена и рассмотрена вырожденная гипергеометрическая функция от двух переменных, через которую выражается решение исследуемого интегрального уравнения Вольтерра первого рода. Найденная формула обращения применена к нахождению некоторых соотношений между искомым решением и его производной краевой задачи для гиперболического уравнения с двумя линиями вырождения и со спектральным параметром. Библиогр.: 16 назв.

**Ключевые слова:** интегральное уравнение Вольтерра первого рода; формула обращения; вырожденная гипергеометрическая функция от двух переменных; спектральный параметр.

**Постановка проблемы и анализ литературы.** Многочисленные приложения теории интегральных уравнений можно найти в теории упругости, в теории пластичности, гидродинамике, в теории массо- и теплопереноса, теории управления, химических технологиях, биомеханике, теории массового обслуживания, экономике и медицине. Часто изучение задач моделирования фильтрации жидкости в пористых средах сводится к рассмотрению интегральных уравнений Вольтера первого рода со специальными функциями в ядрах [1, 2]. Из общей теории известно, что если при решении задач прикладной математики появляется интегральное уравнение Вольтера первого рода, то необходимо найти формулу обращения, выражающую решение интегрального уравнения в явном виде. Подробное изложение теории интегральных уравнений первого рода и библиографию можно найти в [3, 4]. Опуская большой список литературы, в которой изучены различные интегральные уравнения первого рода, отметим работы, наиболее близко примыкающие к настоящему сообщению.

В задаче Коши для гиперболического уравнения с двумя линиями вырождения и со спектральным параметром

$$(-y)^m U_{xx} - x^n U_{yy} + \mu x^n (-y)^m U = 0, \ x > 0, \ y < 0, \qquad (1)$$

где *m*, *n* и μ – действительные числа, при нахождении соотношения между искомым решением и его производной на линии вырождения, т.е. соотношения между $\tau(x) \equiv U(x,+0)$ и $\nu(x) \equiv U'_y(x,+0)$, приходим к следующему интегральному уравнению:

$$N_{0x}^\lambda [\nu(x)] \equiv \int_0^x \left(\frac{t}{x}\right)^\alpha (x-t)^{-2\beta} \times$$
$$\times \Xi_2\left(\alpha,\ 1-\alpha;\ 1-\beta;\ -\frac{(x-t)^2}{4xt}, \lambda(x-t)^2\right) \nu(t) dt, \qquad (2)$$

где $\Xi_2(a,b;d;u,w) = \sum_{m,n=0}^\infty \frac{(a)_m (b)_m}{m! n! (d)_{m+n}} u^m w^n$ – функция Гумберта, $(a)_k$ – символ Похгаммера: $(a)_0 = 1$, $(a)_k = a(a+1)...(a+k-1)$, $k = 1, 2, ...$ а α, β и λ – действительные числа, причем

$$\alpha = \frac{n}{2(n+2)}, \quad \beta = \frac{m}{2(m+2)}, \quad \lambda = \frac{1}{4}\mu.$$

Нетрудно заметить, что при $0 \leq n < 1$ и $0 < m < 1$ параметры α и β принимают положительные значения $0 \leq 2\alpha < 1$ и $0 < 2\beta < 1$, а при $-1 < n \leq 0$ и $-1 < m < 0$ – отрицательные значения $-1 < 2\alpha \leq 0$ и $-1 < 2\beta < 0$, соответственно.

При α = 0 и λ = 0 получим общеизвестное и досконально изученное интегральное уравнение типа Абеля с отрицательным параметром β в виде [3 – 5]

$$\int_0^x (x-t)^{-2\beta} dt = \tau(x).$$

При α = 0 уравнение (2) принимает вид

$$\int_0^x (x-t)^{-2\beta} \overline{J}_{-\beta}[\lambda(x-t)] \nu(t) dt = \tau(x).$$

Это уравнение называется уравнением Вольтерра первого рода с функциями Бесселя в ядрах и оно исследовано многими авторами [3 – 5].

В случае, когда $\lambda = 0$, формула обращения интегрального уравнения с гипергеометрической функцией в ядре

$$\int_0^x (x-t)^{-2\beta} F\left(\alpha, 1-\alpha; 1-\beta; -\frac{(x-t)^2}{4xt}\right) t^\alpha \nu(t) dt = x^\alpha \tau(x)$$

получена в [6, 7].

Многие задачи для уравнений гиперболического типа первого рода с двумя линиями вырождения и со спектральным параметром сводятся к решению интегрального уравнения (2) с неотрицательными параметрами ($0 \le 2\alpha < 2\beta < 1$, $\lambda$ – произвольное число) [8].

Исследованию интегрального уравнения (2) при отрицательных значениях параметров $\alpha$ и $\beta$ посвящены сравнительно мало работ. Когда $-1 < 2\beta < 2\alpha \le 0$ и $\lambda = 0$, отметим лишь работу [9].

При $-1 < n \le 0$ и $-1 < m < 0$ для уравнения (1) прямая $y = 0$ параболического вырождения является особой характеристикой – огибающей обоих семейств характеристик. В зависимости от степеней вырождения $m$ и $n$ предельные значения $\tau(x)$ и $\nu(x)$ могут иметь особенности. Чтобы обеспечить необходимую гладкость решения $U(x, y)$ вне линии характеристического вырождения, необходимо требовать повышенную гладкость функций $\tau(x)$ и $\nu(x)$. С целью ослабить это требование в [10] дано определение и изучены свойства так называемого класса $R_2^\lambda(\alpha, \beta)$ обобщенных решений уравнения (1), который при $\alpha = 0$ и $\lambda = 0$ совпадает с классом $R_2$, введенным и изученным И.Л. Каролем [11]. Кроме того, в [10] на основе известной формулы классического решения задачи Коши [12] для уравнения (1) получен явный вид обобщенного решения этой же задачи в нововведенном классе.

К такому направлению исследований примыкают работы [13, 14].

**Целью** данного исследования является решение интегрального уравнения (2) при $-1 < 2\beta < 2\alpha \le 0$ и при любых значениях $\lambda$ и применение полученную формулу обращения к нахождению некоторых соотношений между $\tau(x)$ и $\nu(x)$.

**Вырожденные гипергеометрические функции от двух переменных.** Нет необходимости говорить о важности свойств гипергеометрических функций. Любой исследователь, имеющий дело с практическими применениями дифференциальных или интегральных уравнений с ними встречается. Решение самых разных задач, относящихся к теплопроводности и динамике, электромагнитным колебаниям и аэродинамике, квантовой механике и теории потенциалов, приводит к изучению гипергеометрических функций.

Разнообразие задач, приводящих к гипергеометрическим функциям, вызвало быстрый рост их числа. Например, в монографии [15] определены области сходимости гипергеометрических функций от трех переменных второго порядка.

В настоящем сообщении мы имеем дело с вырожденными гипергеометрическими функциями от двух переменных. В [15] они определены следующим образом

$$F_{l;m;n}^{p:q;k}\left[\begin{array}{ccc}(a_p): & (b_q); & (c_k); \\ (\alpha_l): & (\beta_m); & (\gamma_n); \end{array} x,y\right] = \sum_{r,s=0}^{\infty} \frac{\prod_{j=1}^{p}(a_j)_{r+s} \prod_{j=1}^{q}(b_j)_r \prod_{j=1}^{k}(c_j)_s}{\prod_{j=1}^{l}(\alpha_j)_{r+s} \prod_{j=1}^{m}(\beta_j)_r \prod_{j=1}^{n}(\gamma_j)_s} \frac{x^r}{r!}\frac{y^s}{s!},$$

где $\prod_{j=1}^{i}(d_j)_t = (d_1)_t (d_2)_t \ldots (d_i)_t$. Область сходимости функции $F_{l;m;n}^{p:q;k}$ известна:

$$p+q < l+m+1, \; p+k < l+n+1, \; |x|<\infty, |y|<\infty$$

или

$$p+q > l+m+1, \; p+k > l+n+1, \; |x|<1, |y|<1,$$

и

$$|x|^{\frac{1}{p-1}} + |y|^{\frac{1}{p-1}} < 1, \text{ если } p>l; \; \max\{|x|,|y|\}<1, \text{ если } p<l.$$

В дальнейшем решение интегрального уравнения (2) будет выражаться с помощью функции

$$F_{1:0;1}^{0:2;1}\left[\begin{array}{ccc}-: & b,c; & d; \\ e: & -; & g; \end{array} x,y\right] = \sum_{m,n=0}^{\infty} \frac{(b)_m (c)_m (d)_n}{(e)_{m+n}(g)_n} \frac{x^m}{m!}\frac{y^n}{n!}.$$

Поэтому опишем некоторые свойства этой функции.

Легко видеть, что функция $F_{1:0;1}^{0:2;1}$ является естественным обобщением известной функции Гумберта $\Xi_2$, т.е. имеет место равенство

$$F_{1:0;1}^{0:2;1}\left[\begin{array}{ccc}-: & b,c; & d; \\ e: & -; & d; \end{array} x,y\right] = \Xi_2(b,c;e;x,y).$$

По методу, изложенному в [16], можно выяснить, что функции $\Xi_2$ и $F_{1:0;1}^{0:2;1}$ являются вырожденными гипергеометрическими функциями второго и третьего порядков [16], соответственно. Известно [16], что функция Гумберта $z = \Xi_2(b,c;e;x,y)$ удовлетворяет систему уравнений

$$\begin{cases} x(1-x)z_{xx} + yz_{xy} + [e-(b+c+1)x]z_x - bcz = 0, \\ yz_{yy} + xz_{xy} + ez_y - z = 0. \end{cases}$$

Далее, следуя работе [16], нетрудно установить, что функция

$$z = F_{1:0;1}^{0:2;1}\begin{bmatrix} - : & b,c; & d; \\ e : & - ; & d; \end{bmatrix} x, y$$

является решением следующей системы уравнений

$$\begin{cases} x(1-x)z_{xx} + yz_{xy} + [e-(b+c+1)x]z_x - bcz = 0, \\ y^2 z_{yyy} + xyz_{xyy} + gxz_{xy} + (e+g+1)yz_{yy} + (eg-y)z_y - dz = 0. \end{cases}$$

Взглянув на эти системы уравнений, легко понять, что порядок вырожденной гипергеометрической функции от двух переменных определяется наивысшим порядком частных производных входящих в данную систему.

**Формула обращения интегрального уравнения (2).** Имеет место следующая

**Теорема.** Пусть $\nu(x)$ – непрерывна и интегрируема в интервале $(0,1)$, $\tau(x) \in C[0,1] \bigcup C^1(0,1)$ и $\tau(0) = 0$. Тогда при $-1 < 2\beta < 2\alpha \leq 0$ и при любых $\lambda$ интегральное уравнение (2) обратимо по формуле

$$\nu(x) = T_{0x}^{\alpha,\beta,\lambda}[\tau(x)] \equiv \frac{\sin 2\beta\pi}{2\beta\pi} x^{-2\alpha} \frac{d}{dx} \left\{ x^\alpha \int_0^x t^\alpha (x-t)^{2\beta} \times \right.$$

$$\left. \times F_{1:0;1}^{0:2;1}\begin{bmatrix} - : & -\alpha, 1+\alpha; & \frac{1}{2}+\beta; \\ 1+\beta : & - ; & -\frac{1}{2}+\beta; \end{bmatrix} -\frac{(x-t)^2}{4xt}, \lambda(x-t)^2 \right] \tau'(t) dt \right\}, \quad (3)$$

и наоборот, т.е. справедливы тождества

$$T_{0x}^{\alpha,\beta,\lambda}\left\{N_{0x}^{\alpha,\beta,\lambda}[\nu(x)]\right\} = \nu(x), \quad N_{0x}^{\alpha,\beta,\lambda}\left\{T_{0x}^{\alpha,\beta,\lambda}[\tau(x)]\right\} = \tau(x). \quad (4)$$

*Доказательство теоремы* осуществим для первого из тождеств (4). Из (2) имеем

$$\tau'(t) = -\beta t^{-\alpha-1} \int_0^t (t-s)^{-2\beta} \Xi_2(\alpha, 1-\alpha; -\beta; u, w) s^\alpha \nu(s) ds -$$

$$- 2\beta t^{-\alpha-1} \int_0^t (t-s)^{-2\beta-1} \Xi_2(\alpha, 1-\alpha; -\beta; u, w) s^{\alpha+1} \nu(s) ds -$$

$$-(\alpha+\beta) t^{-\alpha-1} \int_0^t (t-s)^{-2\beta} F_{1:0;1}^{0:2;1}\left[\begin{array}{c} - : \alpha, 1-\alpha;\ 1-\alpha-\beta; \\ 1-\beta : \quad - \ ;\ -\alpha-\beta; \end{array} u, w\right] s^\alpha \nu(s) ds,$$

где $u = -(t-s)^2/(4ts)$, $w = -\lambda(t-s)^2$.

Подставив теперь найденное выражение для $\tau'(t)$ в (3), после ряда несложных преобразований получим

$$\nu(x) = x^{-2\alpha} \frac{d}{dx}\left\{ x^\alpha \int_0^x W(x,d;\lambda) s^\alpha \nu(s) ds \right\}, \qquad (5)$$

где

$$W(x,s;\lambda) = \sum_{k,n=0}^\infty \frac{(-1)^{k+n}}{k!n!}\left(\beta - \frac{1}{2}\right)_k \left(\frac{1}{2} - \beta\right)_n \Omega(k,n;z) \left[\lambda(x-s)^2\right]^{k+n},$$

$$\Omega(k,n;z) = \sum_{p,q=0}^\infty \frac{(-1)^{p+q}(-\alpha)_p(1+\alpha)_p(\alpha)_q(1-\alpha)_q}{p!q!(1+2k+2n+2p+2q)!} \times$$

$$\times \left(\frac{1}{2}+\beta+k\right)_p \left(\frac{1}{2}-\beta+n\right)_q z^{2p+q}(1-z)^{-q} E(k,n;p,q;z),$$

$$E(k,n;p,q;z) = (-\alpha - 2\beta + 2n + q)zF(2+2k+2n+2p+2q;z) +$$

$$+ (1+2k+2n+2p+2q)(1-z)F(1+2k+2n+2p+2q;z),\ z = (x-s)/x.$$

Здесь для краткости принята запись

$$F(d;z) = F(1+2\beta+2k+2p, 1+p+q; d; z),$$

где $F(a,b;d;u)$ – гипергеометрическая функция Гаусса.

Продолжение доказательства теоремы существенно опирается на свойства функции $W(x,s;\lambda)$. Выделим следующее утверждение в виде леммы, доказательство которой приведем после доказательства теоремы.

**Лемма.** При любых $\lambda$ и $0 < s < x < 1$ справедливо тождество

$$W(x,s;\lambda) = (1-z)^{\alpha}. \qquad (6)$$

Теперь в продолжении доказательства теоремы, если учитывать, что $W(x,s;\lambda) = (s/x)^{\alpha}$, то соотношение (5) превращается в тождество. Тем самым доказано первое из тождеств (3). Аналогично доказывается и второе тождество. Теорема доказана.

*Доказательство леммы.* Пусть $\lambda = 0$. Тогда $W(x,s;0) = \Omega(0,0;z)$. Применяя известную формулу

$$(1-z)^{-q} = \sum_{m=0}^{\infty} \frac{(q)_m}{m!} z^m,$$

получим

$$\Omega(0,0;z) = \sum_{p,q,m,l=0}^{\infty} A(p,q,m,l) z^{2p+q+m+l} + \sum_{p,q,m,l=0}^{\infty} B(p,q,m,l) z^{1+2p+q+m+l},$$

где $A$ и $B$ – известные функции от $p, q, m$ и $l$. Нетрудно установить равенство:

$$\sum_{p,q,m,l=0}^{\infty} A(p,q,m,l) z^{2p+q+m+l} = \sum_{p=0}^{\infty} \sum_{q=0}^{[p/2]} \sum_{m=0}^{p-2q} \sum_{l=0}^{p-2q-m} A(q,m,p-2q-k-l,l) z^p,$$

где $[p]$ – целая часть числа $p$. Используя это равенство и свойства символа Похгаммера, затем приведя подобные члены по степеням $z$, после нескольких преобразований получим

$$W(x,s;0) = (1-z)^{\alpha}. \qquad (7)$$

При наличии $\lambda$, т.е. при $\lambda \neq 0$, следует учитывать, что

$$W(x,s;\lambda) = \sum_{k,n=0}^{\infty} (-1)^k \Omega_1(k;z) \left[\lambda(x-s)^2\right]^k,$$

где

$$\Omega_1(k;z) = \sum_{n=0}^{k} \frac{1}{n!(k-n)!} \left(\beta - \frac{1}{2}\right)_{k-n} \left(\frac{1}{2} - \beta\right)_n \Omega(k-n,n;z).$$

Легко видеть, что $\Omega_1(0;z) = \Omega(0,0;z)$, поэтому функцию $\Omega_1(k;z)$ исследуем при $k \geq 1$.

При каждом $k$ повторив рассуждения, проведенные как в случае $\lambda = 0$, будем иметь $\Omega_1(k;z) = 0$, $k = 1, 2, \ldots$. Следовательно,

$$W(x,s;\lambda) \equiv 0, \quad \lambda \neq 0. \qquad (8)$$

Объединив результаты (7) и (8), заключаем, что при любых $\lambda$ справедливо равенство $W(x,s;\lambda) = (1-z)^{\alpha}$. Лемма доказана.

**Применения операторов (2) и (3).** Рассмотрим уравнение (1) в конечной односвязной области $D$, ограниченной характеристиками

$$AC : \xi \equiv \frac{2}{n+2} x^{(n+2)/2} - \frac{2}{m+2}(-y)^{(m+2)/2} = 0,$$

$$BC : \eta \equiv \frac{2}{n+2} x^{(n+2)/2} - \frac{2}{m+2}(-y)^{(m+2)/2} = 1$$

и $AB : y = 0$ уравнения (1) при $y \leq 0$ и $x \geq 0$, где $m$, $n$ и $\mu$ – действительные числа, причем $-1 < m < 0$, $-1 < n \leq 0$.

В характеристических координатах $\xi$ и $\eta$ уравнение (1) переходит в уравнение типа обобщенного уравнения Эйлера-Пуассона-Дарбу

$$\frac{\partial^2 u}{\partial \xi \partial \eta} + \frac{\alpha}{\eta + \xi}\left(\frac{\partial u}{\partial \eta} + \frac{\partial u}{\partial \xi}\right) - \frac{\beta}{\eta - \xi}\left(\frac{\partial u}{\partial \eta} - \frac{\partial u}{\partial \xi}\right) + \lambda u = 0, \qquad (9)$$

а область $D$ преобразуется в область $\Delta$, граница которой состоит из отрезков прямых $PM : \xi = 0$, $QM : \eta = 1$ и $PQ : \eta = \xi$, где

$$\alpha = \frac{n}{2(n+2)}, \quad \beta = \frac{m}{2(m+2)}, \quad -1 < 2\alpha \leq 0, \quad -1 < 2\beta < 0, \quad \lambda = \frac{1}{4}\mu.$$

Решение задачи Коши для уравнения (9) в области $\Delta$ из класса $R_2^{\lambda}(\alpha, \beta)$ с начальными данными

$$u(\xi, \xi) = \tau(\xi), \quad 0 \leq \xi \leq 1, \qquad (10)$$

$$[2(1-2\beta)]^{-2\beta} \lim_{\eta \to \xi}(\eta - \xi)^{2\beta}\left(\frac{\partial u}{\partial \eta} - \frac{\partial u}{\partial \xi}\right) = \nu(\xi), \quad 0 < \xi < 1 \qquad (11)$$

известно [8]:

$$u(\xi, \eta) = \left(\frac{\eta + \xi}{2}\right)^{-\alpha} \int_0^{\xi} (\eta - t)^{-\beta}(\xi - t)^{-\beta} t^{\alpha} \Xi_2[\alpha, 1-\alpha; 1-\beta; \sigma, \rho] T(t) dt +$$

$$+\left(\frac{\eta+\xi}{2}\right)^{-\alpha}\int\limits_{\xi}^{\eta}(\eta-t)^{-\beta}(t-\xi)^{-\beta}t^{\alpha}\Xi_{2}\bigl[\alpha,1-\alpha;1-\beta;\sigma,\rho\bigr]N(t)dt, \qquad (12)$$

где $T(t)$ – непрерывная и интегрируемая в $(0, 1)$ функция,

$$N(t) = \bigl[2\cos\beta\pi\bigr]^{-1}T(t) - \gamma_{2}\nu(x),$$

$$\gamma_{2} = \bigl[2(1-2\beta)\bigr]^{2\beta-1}\Gamma(2-2\beta)\Gamma^{-1}(1-\beta),$$

$$\sigma = \frac{(\eta-t)(t-\xi)}{2t(\eta+\xi)}, \quad \rho = \lambda(\eta-t)(t-\xi).$$

Для определения $T(t)$ воспользуемся вышеизложенной теоремой. В самом деле, положив $\eta = \xi = x$, из (12) получим

$$\tau(x) = x^{-\alpha}\int\limits_{0}^{x}(x-t)^{-2\beta}t^{\alpha}\Xi_{2}\left[\alpha,1-\alpha;1-\beta;-\frac{(x-t)^{2}}{4xt},-\lambda(x-t)^{2}\right]T(t)dt.$$

Отсюда, в силу формулы (3), найдем

$$T(x) = \frac{\sin 2\beta\pi}{2\beta\pi}x^{-2\alpha}\frac{d}{dx}\Bigl\{x^{\alpha}\int\limits_{0}^{x}t^{\alpha}(x-t)^{2\beta}\times$$

$$\times F_{1:0;1}^{0:2;1}\left[\begin{array}{cc} - : & -\alpha,1+\alpha; \quad \frac{1}{2}+\beta; \\ 1+\beta : & - \quad ; \quad -\frac{1}{2}+\beta; \end{array} -\frac{(x-t)^{2}}{4xt},-\lambda(x-t)^{2}\right]\tau'(t)dt\ \Bigr\}.$$

Приведем другое применение формул обращения (2), (3).

Рассмотрим задачу Коши-Гурса для уравнения (9) с условиями (11) и

$$u(0,\eta) = \varphi(\eta), \ 0 \le \eta \le 1.$$

Решение этой задачи из класса $R_{2}^{\lambda}(\alpha,\beta)$ имеет явный вид [14]:

$$u(\xi,\eta) = \left(\frac{\eta+\xi}{2}\right)^{-\alpha}\int\limits_{0}^{\xi}(\eta-t)^{-\beta}(\xi-t)^{-\beta}t^{\alpha}\Xi_{2}\bigl[\alpha,1-\alpha;1-\beta;\sigma,\rho\bigr]\Psi(t)dt +$$

$$+\left(\frac{\eta+\xi}{2}\right)^{-\alpha}\int\limits_{\xi}^{\eta}(\eta-t)^{-\beta}(t-\xi)^{-\beta}t^{\alpha}\Xi_{2}\bigl[\alpha,1-\alpha;1-\beta;\sigma,\rho\bigr]\Phi(t)dt, \qquad (13)$$

где $\Psi(x) = 2\gamma_2 \cos\beta\pi \cdot \nu(x) + \Phi(x)$, $\Phi(x)$ – известная функция.

Формула (13) играет важную роль при изучении задач для уравнений смешанного типа, так как из нее при $\eta = \xi = x$ легко вывести основное функциональное соотношение между $\tau(x)$ и $\nu(x)$ на линии вырождения, принесенное из гиперболической части смешанной области.

Действительно, предположим в задаче Коши-Гурса $\varphi(x) = 0$ и в (13) положим $\eta = \xi = x$. Тогда получаем интегральное уравнение вида (2). Теперь воспользовавшись теоремой, решение последнего уравнения находим в виде (3).

**Выводы.** Таким образом, в результате исследований получена формула обращения (3) интегрального уравнения (2) при $-1 < 2\beta < 2\alpha \leq 0$ и произвольных значениях параметра $\lambda$. Формула (3) является очень важным инструментом при исследовании локальных и нелокальных краевых задач для уравнений смешанных параболо-гиперболического и эллиптико-гиперболического типов второго рода со спектральным параметром.

Ehrgashev Tuhtasin, Cand. Phy.-Math. Sci., Associate professor
Tashkent Institute of Irrigation and Agricultural Mechanization Engineers
Uzbekistan, Tashkent, 100000, Kari-Niyazi st., 39
Tel.: +99894 673 1869, +99871 2230524, E-mail: ertuhtasin@mail.ru